\documentclass[12pt]{amsart}

\usepackage{mathrsfs}

\usepackage{amsfonts, amssymb,amsmath, amsthm, amsxtra, latexsym, amscd}
\usepackage[latin1]{inputenc} 

\usepackage[english]{babel}

\usepackage{amsmath}

\usepackage{amsthm,amssymb,latexsym}

\usepackage{t1enc}

\usepackage{epic}

\usepackage{eepic}

\usepackage{xypic}

\usepackage{amsfonts, amssymb,amsmath, amsthm, amsxtra, latexsym, amscd}

\usepackage{graphics}

\begin{document}

\title[Tube domains and restrictions of minimal representations]{Tube domains and restrictions of minimal representation}
\author{Henrik Seppänen}
\address{Department of Mathematics, Chalmers University of Technology and G\"oteborg University, G\"oteborg, Sweden}
\email{henriks@math.chalmers.se}
\keywords{Lie groups, unitary representations, branching law, real bounded symmetric domains}
\subjclass{22E45, 32M15, 33C45, 43A85}

\newtheorem{prop}{Proposition}

\newtheorem{lemma}[prop]{Lemma}

\newtheorem{cor}[prop]{Corollary}

\newtheorem{thm}[prop]{Theorem}

\newtheorem{con}[prop]{Conjecture}

\theoremstyle{definition}

\newtheorem{defin}[prop]{Definition}

\theoremstyle{remark} \newtheorem*{remark}{Remark}

\bibliographystyle{amsplain}

\thispagestyle{empty}






\maketitle

\begin{abstract}
In this paper we study the restrictions of the minimal representation in the analytic continuation of the scalar holomorphic discrete
series from $Sp(n,\mathbb{R})$ to $GL(n,\mathbb{R})$, and from $SU(n,n)$ to $GL(n,\mathbb{C})$ respectively.
We work with the realisations of the representation spaces as $L^2$-spaces on the boundary orbits of rank one of the 
corresponding cones, and give explicit integral operators that play the role of the intertwining operators for the decomposition.
We prove inversion formulas for dense subspaces and use them to prove the Plancherel theorem for the respective decomposition.
The Plancherel measure turns out to be absolutely continuous with respect to the Lebesgue measure in both cases.
\end{abstract}

\section{Introduction}
The unitary representations obtained by continuation of the scalar holomorphic discrete series of a hermitian Lie group, $G$,
were classified by 
Wallach in \cite{wall}, and independently by Rossi and Vergne (\cite{Rossi-V}). The classification amounts
to membership in the \emph{Wallach set} for the linear functionals on the compact Cartan subalgebra that extend the 
family of weights parametrising the weighted Bergman spaces on the symmetric space $G/K$. 

These unitary representations can all be realised on Hilbert spaces of holomorphic functions on the 
corresponding bounded symmetric domain $\mathscr{D} \cong G/K$.
However, in this model the unitary structure cannot be described in a uniform way even though the corresponding reproducing kernels can.
In any case, the restriction to any totally real submanifold defines an injective mapping. Therefore it is natural to
consider an antiholomorphic involution $\tau: \mathscr{D} \rightarrow \mathscr{D}$ that lifts to an involutive automorphism 
(which we also denote by $\tau$) of the group $G$. Letting $H=G^{\tau}$ denote the fixed point group, and $L=K \cap H$, the
space $\mathscr{X}:=H/L$ is a totally real submanifold.
The decomposition of the restriction to $H$ of the unitary representations obtained by analytic continuation, or more generally, 
the restriction of holomorphic representations to symmetric subgroups, has lately been an area of intensive research.
Among those who have studied this problem we find, for example, Davidson, Ólafsson, and Zhang (\cite{doz}), van Dijk and Pevzner 
(\cite{dijkpevz}), Zhang
(\cite{zhtams}, \cite{zhtp}, \cite{zhSB}), and the author (\cite{papper1}, \cite{papper3}).
However, there does not yet seem to be any uniform way of dealing with this problem. For regular parameter, the Segal-Bargmann
transform provides a unitary equivalence between the restriction to the group $H$ and the left regular representation of $H$
on the space $L^2(H/L)$ and in this case the decomposition is determined explicitly by the Helgason Fourier transform for
the symmetric space $H/L$. Otherwise, the results obtained so far depend on particular features of the special cases.
In \cite{zhtp}, Zhang decomposes the restriction to the diagonal subgroup of the tensor product of a minimal representation and its dual
by finding the spectral decomposition for the Casimir operator.
The same method is used in \cite{papper3}, where the author determines the restriction of the minimal representation for $SU(n,m)$
to the subgroup $SO(n,m)$. It should be noted that this approach identifies the representations occurring in the decomposition
and determines the Plancherel measure explicitly, but it does not provide an intertwining operator.
In \cite{papper1}, the author determines the restriction from $SO(2,n)$ to $SO(1,n)$ for general parameter in the Wallach set
and gives an intertwining operator. This was possible thanks to an explicit power series expansion for the spherical functions
on the group $SO(1,n)$.

In this paper we consider the minimal representations for the groups $Sp(n, \mathbb{R})$ and $SU(n,n)$ and restrict
to the automorphism groups for the cones associated with the respective tube domains. We use the model
in \cite{Rossi-V} that realises the representations as $L^2$-spaces on the orbits of rank one elements on the boundaries
of the respective cones. The intertwining operators are given explicitly as integral transforms (or, rather, as 
analytic continuations of operators defined as integral transforms).
The two cases we deal with are identical in principle. However, the proofs are rather technical when it comes to parameters, so 
we chose to avoid a uniform presentation to increase the readability. Instead we present two parallel cases where
the solutions follow the same guideline.

The paper is organised as follows. Section 2 contains preliminaries for the two cases separately. In section 3 we describe
the constituents in the decomposition for the restriction from $Sp(n,\mathbb{R})$, construct an intertwining operator and
prove the Plancherel theorem. Section 4 is the analogue of section 3 for the group $SU(n,n)$.

{\bf Acknowledgement:} The author is grateful to his advisor Professor Genkai Zhang for illuminating discussions on the
topic of this paper.
  
\section{Preliminaries}

\subsection{\emph{Type $II_n$}}
Let $V$ be the real vector space of symmetric $n \times n$-matrices. The complexification, 
$V^{\mathbb{C}}=V \oplus iV$, of $V$ 
consists of all complex symmetric $n \times n$ matrices. 
Consider the bounded symmetric domain 

\begin{eqnarray}
\mathscr{D}=\{Z \in V^{\mathbb{C}} | I-Z^*Z>0\}.
\end{eqnarray}
The group 

\begin{eqnarray*}
G=Sp\,(n, \mathbb{R})=\left\{g \in SU(n,n) | g^t 
\left( \begin{array}{cc}
0 & I\\
-I & 0 \\
\end{array}\right)g=
\left( \begin{array}{cc}
0 & I\\
-I & 0 \\
\end{array}\right)
 \right\}
\end{eqnarray*}
acts transitively on $\mathscr{D}$ by

\begin{eqnarray}
Z \mapsto (AZ+B)(CZ+D)^{-1},
\end{eqnarray}
where 
\begin{eqnarray}
g= \left(
\begin{array}{cc}
A & B\\
C & D\\
\end{array}\right)
\end{eqnarray}
consists of the $n \times n$ blocks $A, B, C$ and $D$. 
The isotropy group of $0$ is 

\begin{eqnarray}
K=\left\{ \left(
\begin{array}{cc}
A & 0 \\
0 & \overline{A}\\
\end{array}\right)| A \in U(n) \right\},
\end{eqnarray}
and hence

\begin{eqnarray}
\mathscr{D} \cong G/K.
\end{eqnarray}
Let 
$$\Omega=\{X \in V| X >0\}.$$ 
Then $\Omega$ is a symmetric cone in $V$ with automorphism group 
$GL(n,\mathbb{R})$ acting as $$X \stackrel{g}{\mapsto} gXg^t.$$
In fact $$\Omega \cong GL(n, \mathbb{R})/O(n).$$
Moreover, the boundary of $\Omega$ is partitioned into $n$ orbits under $GL(n, \mathbb{R})$,
$$\partial \Omega = \cup_{i=1}^n \Omega^{(i)},$$
where $\Omega^{(i)}$ is the set of positive semidefinite matrices of rank $i$. Each orbit carries a 
\emph{quasi-invariant} measure, $\mu_i$, transforming in the fashion
$$g^*\mu_i=| \det g|^i \mu_i$$ under the action of $GL(n, \mathbb{R})$.

The Cayley transform $$c(Z)=(I-Z)(I+Z)^{-1}$$ maps $\mathscr{D}$ biholomorphically onto the tube domain
\begin{eqnarray*}
T_{\Omega}:=\{Z=U+iV \in V^{\mathbb{C}}| U \in \Omega \}.
\end{eqnarray*}
Let $\tau$ denote the conjugation with respect to $V$, i.e., $$\tau(u+iv)=u-iv.$$
The set, $\mathscr{X}$ of fixed points of $\tau$ in $\mathscr{D}$, $V \cap \mathscr{D}$ is a totally real and totally geodesic 
real submanifold of $\mathscr{D}$, and the Cayley transform restricts to a diffeomorphism
$$\mathscr{X} \cong \Omega.$$ In particular, $\mathscr{X}$ is a homogeneous space $$\mathscr{X} \cong H/L,$$ where 
$H \cong GL(n, \mathbb{R})$ and $L \cong O(n)$.
Consider now the isomorphism $$\Psi : GL(n, \mathbb{R}) \rightarrow \mathbb{R}^* \times SL(n,\mathbb{R})$$ given by
$$g \stackrel{\Psi}{\mapsto} (\mbox{det}(g),\mbox{det}(g)^{-1/n}g)$$ with inverse $\Psi^{-1}$ given by 
$$(\lambda, h) \stackrel{\Psi^{-1}}{\mapsto} \lambda^{1/n}h.$$
The differential of $\Psi$ at the identity element gives an isomorphism of Lie algebras 
$$\mathfrak{gl}(n, \mathbb{R})=\mathbb{R} \oplus \mathfrak{sl}(n,\mathbb{R}).$$
We will denote the $SL(n, \mathbb{R})$-factor in $H$ by $H'$, and correspondingly we let $\mathfrak{h'}$ denote  
$\mathfrak{sl}(n,\mathbb{R})$.

The minimal representation in the analytic continuation of the scalar holomorphic discrete series of $G$ 
can be defined as a Hilbert space of functions, $\mathscr{H}_{1/2}$, on $T_{\Omega}$. This space has the
reproducing kernel
\begin{eqnarray}
K_{1/2}(z,w):=\det(z-w^*)^{-1/2}.
\end{eqnarray}
By (an analytic continuation of) a Laplace transform, it is unitarily and $G$-equivariantly equivalent to the
the Hilbert space $L^2(\Omega^{(1)}, \mu_1)$ (cf. \cite{Rossi-V}). Another proof of this can be found in \cite{FKsymm}. 
We will now give an even more explicit model 
for this representation space.
Consider therefore the mapping 
$$\eta : \mathbb{R}^n \setminus \{0\} \rightarrow \Omega^{(1)},$$
defined by $$\eta(x)=xx^t.$$ Here we identify $\mathbb{R}^n$ with the space of all $n \times 1$ real matrices. 
It is straigthforward to check that $\eta$ is surjective and that $$\eta(x)=\eta(y) \Leftrightarrow 
x=\pm y,$$ and hence we have a bijection $$\Omega^{(1)} \cong (\mathbb{R}^n \setminus \{0\})/ \pm 1,$$ where the right hand side denotes the 
set of orbits under the linear action of the two-element group generated by the endomorphism $-I$.
Moreover, the action of $GL(n, \mathbb{R})$ is covered by the linear action on $\mathbb{R}^n \setminus \{0\}$ so that we have the following
commuting diagram.

\begin{eqnarray*}
\begin{CD}
 \mathbb{R}^n \setminus \{0\} @>x \mapsto gx>> \mathbb{R}^n \setminus \{0\}  \\
\eta@VVV      @VV\eta V\\
\Omega^{(1)} @>X \mapsto gXg^t>> \Omega^{(1)}
\end{CD}
\end{eqnarray*}

The measure $\mu_1$ is the pushforward under $\eta$ of the Lebesgue measure on $\mathbb{R}^n \setminus \{0\}$.
Hence, the minimal representation can be realised in the Hilbert space of even square-integrable functions on 
$\mathbb{R}^n \setminus \{0\}$. In this picture we have the formula 
\begin{eqnarray}
f \stackrel{h}{\mapsto} \det h \,(f \circ h^t)
\end{eqnarray}
for the group action on 
functions. Obviously this Hilbert space contains $L$-invariant functions, i.e., the representation is spherical. 
Therefore, by \cite{papper1}, there exists a direct integral decomposition 
$$L^2((\mathbb{R}^n \setminus \{0\})/\pm 1) \cong \int_{\Lambda}\mathcal{H}_{\lambda}d\nu,$$ where
$\Lambda$ is some parameter set, the $\mathcal{H}_{\lambda}$ are canonical representation spaces for irreducible spherical 
unitary representations of $H'$, and $\nu$ is some positive measure on $\Lambda$, called the \emph{Plancherel measure} for
the minimal representation. 

\subsection{\emph{Type $I_{nn}$}}
Let $V$ be the real vector space of Hermitian $n \times n$-matrices. The complexification, 
$V^{\mathbb{C}}=V \oplus iV$, of $V$ 
consists of all complex $n \times n$ matrices. 
Consider the bounded symmetric domain
\begin{eqnarray*} 
\mathscr{D}=\{Z \in V^{\mathbb{C}} | I-Z^*Z>0\}.
\end{eqnarray*}
The group $$G=SU(n,n)$$ acts on $\mathscr{D}$ by 
$$Z \mapsto (AZ+B)(CZ+D)^{-1},$$
where $$g= \left(
\begin{array}{cc}
A & B\\
C & D\\
\end{array}\right)$$ consists of the $n \times n$ blocks $A, B, C$ and $D$. 
We have the description
\begin{eqnarray*}
\mathscr{D} \cong G/K.
\end{eqnarray*}
of $\mathscr{D}$ as a homogeneous space, where
\begin{eqnarray*}
K=S(U(n) \times U(n)
= \left\{ \left(
\begin{array}{cc}
A & 0 \\
0 & D\\
\end{array}\right)| A , D \in U(n), \mbox{det}(A)\mbox{det}(D)=1 \right\}.
\end{eqnarray*}
The symmetric cone
\begin{eqnarray*}
\Omega=\{X \in V| X >0\}
\end{eqnarray*}
in $V$ has automorphism group 
$GL(n,\mathbb{C})$ acting as 
\begin{eqnarray*}
X \stackrel{g}{\mapsto} gXg^*,
\end{eqnarray*}
and
\begin{eqnarray*}
\Omega \cong GL(n, \mathbb{C})/U(n).
\end{eqnarray*}
The boundary of $\Omega$ partitions into $GL(n, \mathbb{C})$ orbits as
$$\partial \Omega = \cup_{i=1}^n \Omega^{(i)},$$
where $\Omega^{(i)}$ is the set of positive semidefinite matrices of rank $i$. Each orbit carries a 
\emph{quasi-invariant} measure, $\mu_i$, transforming in the fashion
$$g^{*} \mu_i =| \det g|^{2i} \mu_i$$ under the action of $GL(n, \mathbb{C})$. 
The Cayley transform $$c(Z)=(I-Z)(I+Z)^{-1}$$ maps $\mathscr{D}$ biholomorphically onto the tube domain
\begin{eqnarray*}
T_{\Omega}:=\{Z=U+iV \in V^{\mathbb{C}}| U \in \Omega \}.
\end{eqnarray*}
Let $\tau$ denote the conjugation with respect to $V$, i.e., $$\tau(u+iv)=u-iv.$$
The set, $\mathscr{X}$ of fixed points of $\tau$ in $\mathscr{D}$, $V \cap \mathscr{D}$ is a totally 
real submanifold of $\mathscr{D}$ and the Cayley transform restricts to a diffeomorphism
$$\mathscr{X} \cong \Omega.$$ In particular, $\mathscr{X}$ is a homogeneous space $$\mathscr{X} \cong H/L,$$ where 
$H \cong GL(n, \mathbb{C})$ and $L \cong U(n)$.

The minimal representation in the analytic continuation of the scalar holomorphic discrete series of $G$ 
is a Hilbert space, $\mathscr{H}_1$, of functions on $T_{\Omega}$ with reproducing kernel
\begin{eqnarray*}
K_1(z,w):=\det(z-w^*)^{-1}.
\end{eqnarray*}
Another realisation is furnished by
the Hilbert space $L^2(\Omega^{(1)}, \mu_1)$ (cf. \cite{Rossi-V}, \cite{FKsymm}). 

To give an explicit realisation, we now consider the mapping 
\begin{eqnarray*}
\eta : \mathbb{C}^n \setminus \{0\} \rightarrow \Omega^{(1)}
\end{eqnarray*}
defined by $$\eta(z)=zz^*.$$ Here we identify $\mathbb{C}^n$ with the space of all $n \times 1$ complex matrices. 
It is straigthforward to check that $\eta$ is surjective and that $$\eta(z)=\eta(w) \Leftrightarrow 
z=e^{i\theta }w,$$ for some real $\theta$, and hence we have a bijection 
$$\Omega^{(1)} \cong \mathbb{C}^n/{U(1)}.$$
The action of $GL(n, \mathbb{C})$ is covered by the linear action on $\mathbb{C}^n$ so that we have the  
following commuting diagram.

\begin{eqnarray*}
\begin{CD}
 \mathbb{C}^n \setminus \{0\} @>z \mapsto gz>> \mathbb{C}^n \setminus\{0\} \\
\eta@VVV      @VV\eta V\\
\Omega^{(1)} @>Z \mapsto gZg^*>> \Omega^{(1)}
\end{CD}
\end{eqnarray*}

The measure $\mu_1$ is the pushforward under $\eta$ of the Lebesgue measure on $\mathbb{C}^n \setminus \{0\}$.
Hence, the minimal representation can be realised in the Hilbert space of $U(1)$-invariant square-integrable functions on 
$\mathbb{C}^n$. In this picture, we have the formula 
\begin{eqnarray}
f \stackrel{h}{\mapsto} \mbox{det}_{\mathbb{R}}h^*\,(f \circ h^*),
\end{eqnarray}
for the group action on 
functions, where the subscript on the determinant means the determinant of $h$ as an $\mathbb{R}$-linear operator on $\mathbb{R}^{2n}$.

\section{The branching rule: type $II_n$}

\subsection{\emph{Some parabolically induced representations}}
In the following, we will consider some parabolically 
induced representations of
$H'=SL(n,\mathbb{R})$.

Let $\mathfrak{a}_0=\mathbb{R}e$, where 

\begin{eqnarray*}
e=
\left(\begin{array}{cc}
n-1 & 0 \\
0 & -I_{n-1}
\end{array} \right),
\end{eqnarray*}
where $I_{n-1}$ denotes the identity matrix of size $(n-1)\times(n-1)$.
The maximal parabolic subalgebra, $\mathfrak{q_0}$, determined by $\mathfrak{a}_0$ has a decomposition
$$\mathfrak{q_0}=\overline{\mathfrak{n}_0} \oplus \mathfrak{m}_0 \oplus \mathfrak{a}_0 \oplus \mathfrak{n}_0,$$
where
\begin{eqnarray*}
\mathfrak{n_0}&=&
\left\{ \left(\begin{array}{cccc}
0 & x_1 & \cdots & x_{n-1} \\
0 & 0 & \cdots & 0\\
\vdots & \vdots & \ddots & 0\\
0 & 0 & \cdots & 0
\end{array} \right)
| x_1, \ldots ,x_{n-1} \in \mathbb{R} \right \},
\\
\mathfrak{m_0}&=&
\left\{ \left(\begin{array}{cc}
0 & 0\\
0 & M 
\end{array} \right)
| M \in \mathfrak{sl}(n-1,\mathbb{R}) \right \},
\\
\overline{\mathfrak{n_0}}&=&
\left\{ \left(\begin{array}{cccc}
0 & 0 & \cdots & 0 \\
x_1 & 0 & \cdots & 0\\
\vdots & \vdots & \ddots & 0\\
x_{n-1} & 0 & \cdots & 0
\end{array} \right)
| x_1, \ldots ,x_{n-1} \in \mathbb{R} \right \}.
\end{eqnarray*}
Here the subspace $\mathfrak{m}_0$ is defined by the property
$$Z_{\mathfrak{h'}}(\mathfrak{a}_0)=\mathfrak{a}_0 \oplus \mathfrak{m}_0,$$
and
\begin{eqnarray*}
\mathfrak{n}_0&=&\left \{ X \in \mathfrak{h'} |\quad [H,X]=\alpha(H)X, \quad \forall H \in \mathfrak{a}_0 \right\},\\
\overline{\mathfrak{n}_0}&=&\left \{ X \in \mathfrak{h'} |\quad [H,X]=-\alpha(H)X, \quad \forall H \in \mathfrak{a}_0 \right\}
\end{eqnarray*}
are the generalised root spaces, 
where the root $\alpha \in \mathfrak{a_0}^*$ is determined by 
$$\alpha(e)=n.$$
We let $\rho_0$ denote the half sum of the positive roots counted with multiplicity, i.e., 
$$\rho_0=\frac{n-1}{2}\alpha.$$
On the group level we have the corresponding decomposition
$$Q_0=M_0A_0N_0,$$
where
\begin{eqnarray*}
A_0 &=&
\left \{ \left ( \begin{array}{cc}
e^s & 0 \\
0 & qI_{n-1}\\
\end{array} \right) | s,q\in \mathbb{R}, e^sq^{n-1}=1\right\}\\
M_0 &=& 
\left\{ \left(\begin{array}{cc}
1 & 0\\
0 & M 
\end{array} \right)
| M \in GL(n-1, \mathbb{R}) \right \},\\
N_0 &=&
\left\{ \left( \begin{array}{ccccc}
1 & x_1 & x_2 & \cdots & x_{n-1}\\
0 & 1 & 0 & \cdots & 0\\
0 & 0 & 1 & \cdots & 0\\
\vdots & \vdots & 0 & \ddots & 0\\
0 & 0 & 0 & \cdots & 1
\end{array}
\right)| x_1, \ldots ,x_{n-1} \in \mathbb{R} \right\}
\end{eqnarray*}

Consider now the representation $1 \otimes \exp i\lambda \otimes 1$ of the group 
\begin{eqnarray*}
Q_0=M_0A_0N_0.
\end{eqnarray*}
The induced representation 
\begin{eqnarray}
\tau_{\lambda}:=\mbox{Ind}_{Q_0}^{H'}(1 \otimes \exp (i\lambda+\rho_0) \otimes 1)
\end{eqnarray}
has a noncompact realisation
in the Hilbert space $L^2(\overline{N_0}, d\overline{n})$ (cf. \cite{knapp1}).    
We have 
\begin{eqnarray}
\pi_{\lambda}(h) f(\overline{n})=
e^{-(i\lambda+\rho_0)(\log a_0(h^{-1}\overline{n}))}f(\overline{n_0}(h^{-1}\overline{n})),
\end{eqnarray}
The decomposition $\mathfrak{h'}=\overline{\mathfrak{n_0}} \oplus \mathfrak{m_0} \oplus \mathfrak{a_0} \oplus {n_0}$
gives a corresponding decomposition $$H'\stackrel{.}{=}\overline{N_0}M_0A_0N_0,$$ by which we mean that the equality holds
outside a set of strictly lower dimension. The factorisation of a group element with respect to this decomposition is not unique, but
the $A_0$-component is. For $h \in H'$, we let $a_0(h)$ denote this component.
The mapping $\exp: \mathfrak{a}_0 \rightarrow A_0$ is a diffeomorphic homomorphism of abelian groups. We let $\log: A_0 \rightarrow 
\mathfrak{a}_0$ denote its inverse.
The representation $\tau_{\lambda}$ is then given by
\begin{eqnarray}
\tau_{\lambda}(h) f(\overline{n})=
e^{-(i\lambda+\rho_0)(\log a_0(h^{-1}\overline{n}))}f(\overline{n_0}(h^{-1}\overline{n})).
\end{eqnarray}
For arbitrary $h \in H', \overline{n} \in \overline{N_0}$, a factorisation of $h\overline{n}$ can be given by
\begin{eqnarray*}
&&\left(\begin{array}{cc}
a & b\\
c & d \end{array} \right)
\left( \begin{array}{cc}
1 & 0 \\
x & I_{n-1}\\
\end{array} \right)\\
&&=\left(\begin{array}{cc}
1 & 0 \\
\frac{c+dx}{a+bx} & I_{n-1} 
\end{array} \right)
\left( \begin{array}{cc}
\frac{a+bx}{|a+bx|} & 0\\
0 & |a+bx|^{1/n-1}(d-\left(\frac{c+dx}{a+bx}\right)b)
\end{array} \right)\\
&&\times \left(\begin{array}{cc}
|a+bx| & 0\\
0 & |a+bx|^{-1/n-1}I_{n-1} 
\end{array} \right)
\left( \begin{array}{cc}
1 & \frac{b^t}{a+bx}\\
0 & I_{n-1}
\end{array} \right). 
\end{eqnarray*}

In view of this, and identifying $L^2(\overline{N_0},d\overline{n})$ with $L^2(\mathbb{R}^{n-1},dx)$,
we obtain the following explicit formula for the induced representation:
\begin{eqnarray*}
(\tau_{\lambda}(h) f)(x)=|ax+b|^{-(i\lambda +n/2)}f \left(\frac{c+dx}{a+bx}\right ),
\end{eqnarray*} 
where $h=\left(\begin{array}{cc}
a & b\\
c & d \end{array} \right)$.

\subsection{\emph{An intertwining operator}}
Let $m_{\lambda}$ be the one-dimensional representation
\begin{eqnarray*}
m_{\lambda}(c)z:=|c|^{i\lambda+\rho_0}z
\end{eqnarray*}
of $\mathbb{R}^*$.
Recalling the isomorphism  
\begin{eqnarray*}
\Psi : GL(n, \mathbb{R}) \rightarrow \mathbb{R}^* \times SL(n,\mathbb{R})
\end{eqnarray*}
from the previous section, 
we can now form 
the representation 
\begin{eqnarray}
\pi_{\lambda}:=m_{\lambda} \otimes \tau_{\lambda}
\end{eqnarray}
of $GL(n, \mathbb{R})$. We let $\mathcal{H}_{\lambda}$ denote the associated representation space.

We now consider an operator, $T$, mapping a $C_0^{\infty}((\mathbb{R}^n \setminus \{0\})/\pm1)$-function, $f$, to a function 
$Tf:\mathbb{C} \times \mathbb{R}^{n-1} \rightarrow \mathbb{R}$ which is meromorphic in the first variable.
The function 
\begin{eqnarray}
Tf(z,\eta):=\int_{\mathbb{R}^n}f(x)|\langle x, (1, \eta)\rangle|^{-(iz+n/2)}dx \label{T-integral}
\end{eqnarray}
is well defined as a function of $z$ and $\eta$ when $\mbox{Im} z-n/2>-1$. For such $z$, a change of variables, followed
by an integration by parts, yields 
\begin{eqnarray*}
&&\int_{\mathbb{R}^n}f(x)|\langle x, (1, \eta)\rangle|^{-(iz+n/2)}dx\\
&=&\frac{1}{(1+|\eta|^2)^{iz+n/2}(-(iz+n/2)+1)} \\
&&\times \int_{\mathbb{R}^{n-1}}\int_{y_1<0}\frac{\partial (f \circ g)(y)}{\partial y_1}|y_1|^{-(iz+n/2)+1}dy_1dy_2 \ldots dy_n\\
&&-\frac{1}{(1+|\eta|^2)^{iz+n/2}(-(iz+n/2)+1)}  \\
&&\times \int_{\mathbb{R}^{n-1}}\int_{y_1>0}\frac{\partial (f \circ g)(y)}{\partial y_1}|y_1|^{-(iz+n/2)+1}dy_1dy_2 \ldots dy_n,
\end{eqnarray*}
where $g$ is some orthogonal transformation such that
\begin{eqnarray*}
(1,\eta)=g(|(1,\eta)|e_1).
\end{eqnarray*}
 By repeated integration by parts we get the identity
\begin{eqnarray}
&&Tf(z,\eta) \label{meromorf} \\
&=&\frac{1}{(1+|\eta|^2)^{iz+n/2}\Pi_{j=1}^k(-(iz+n/2)+j)} \nonumber \\
&& \times \int_{\mathbb{R}^{n-1}}\int_{y_1<0}\frac{\partial^k (f\circ g)(y)}{\partial y_1^k}|y_1|^{-(iz+n/2)+k}dy_1dy_2 \ldots dy_n\nonumber \\
&&+(-1)^k\frac{1}{(1+|\eta|^2)^{iz+n/2}\Pi_{j=1}^k(-(iz+n/2)+j)}\nonumber \\
&&\times \int_{\mathbb{R}^{n-1}}\int_{y_1>0}\frac{\partial^k (f \circ g)(y)}{\partial y_1^k}|y_1|^{-(iz+n/2)+k}dy_1dy_2 \ldots dy_n. \nonumber
\end{eqnarray}
Therefore, $Tf(z,\eta)$ can be continued to a meromorphic function with poles at $z=i(n/2-j)$, for $j=1,2, \ldots$.
In particular, $Tf(\lambda,\eta)$ is a well defined real analytic function for real $\lambda$ (except possibly for $\lambda=0$).
We write $T_{\lambda}f(\eta)$ for $Tf(\lambda,\eta)$. 

\begin{prop} \label{L2}
For $f \in C_0^{\infty}((\mathbb{R}^n \setminus \{0\}) / \pm 1)$ and $\lambda \in \mathbb{R}$, the functions $T_{\lambda}f$ is
in $L^2(\mathbb{R}^{n-1})$.
\end{prop}

\begin{proof}
For $\lambda \in \mathbb{R}$, choose the natural number $k>n/2$ in \eqref{meromorf}. The functions
\begin{eqnarray}
\frac{\partial^k}{\partial y_1^k}(f \circ h), \qquad h \in O(n) 
\end{eqnarray}
constitute a uniformly bounded family in the supremum-norm.
Hence, we have an estimate
\begin{eqnarray}
|T_{\lambda}f(\eta)| \leq C(\lambda)(1+|\eta|^2)^{-(i\lambda+n/2)},
\end{eqnarray}
and this proves the claim.
\end{proof}

In what follows, we will state and prove properties for the functions
$T_{\lambda}f$ for arbitrary real $\lambda$ although the proofs will use the defining integral \eqref{T-integral} which 
makes sense only when $\mbox{Im}z>n/2-1$. The idea is then that both sides in the stated
equalities are meromorphic functions, so by the uniqueness theorem for meromorphic functions it suffices 
to perform the calculations when the defining integral makes sense. All integral equalities should therefore be thought of as analytic
continuations of the corresponding equalities when the integrals are convergent.

\begin{prop} \label{H-ekv1}
The operator 
$$T_{\lambda}: C_0^{\infty}(\mathbb{R}^n/ \{\pm 1\}) \rightarrow \mathcal{H}_{\lambda},$$
given by
\begin{eqnarray}
T_{\lambda}f(\eta)=\int_{\mathbb{R}^n}f(x)|\langle x, (1, \eta)\rangle|^{-(i\lambda+n/2)}dx
\end{eqnarray}
is $H$-equivariant.
\end{prop}

\begin{proof}
Take $g \in H$ and write $g=\zeta h$, where $\zeta$ is a diagonal matrix and $h$ has determinant 1. Moreover we write 
$h^{-1}=
\left( \begin{array}{cc}
a & b\\
c & d
\end{array} \right)$.
Then
\begin{eqnarray*}
&&T_{\lambda}(gf)(\eta)\\
&=&\int_{\mathbb{R}^n}f(g^tx)|\langle(x_1,x'),(1,\eta)\rangle|^{-(i\lambda+n/2)}\det g dx\\
&=&\int_{\mathbb{R}^n}f(x)|\langle(x_1,x'),g^{-1}(1,\eta)\rangle|^{-(i\lambda+n/2)}dx\\
&=&|\zeta|^{i\lambda+n/2}\int_{\mathbb{R}^n}f(x)|(a+b\eta)x_1+\langle x',c+d\eta \rangle|^{-(i\lambda+n/2)}dx\\
&=&|\zeta|^{i\lambda+n/2}\int_{\mathbb{R}^n}f(x)|a+b\eta|^{-(i\lambda+n/2)}\\
&& \qquad \qquad \times |\langle(x_1,x'),(1,(c+d\eta)(a+b\eta)^{-1})\rangle|^
{-(i\lambda+n/2)}dx\\
&=&\pi_{\lambda}(g)T_{\lambda}f(\eta).
\end{eqnarray*}
\end{proof}
If $f$ is $L$-invariant, then $T_{\lambda}f$ is an $L$-invariant function in the representation space $\mathcal{H}_{\lambda}$. 
By the Cartan-Helgason theorem (\cite{helg2}), the
subspace of $L$-invariants is at most one dimensional. In fact, it is spanned by the function 
$\eta \mapsto (1+|\eta|^2)^{-(i\lambda+n/2)/2}$. Thus, we can define a function $\tilde{f}$ by 
\begin{eqnarray}
T_{\lambda}f(\eta)=\tilde{f}(\lambda)(1+|\eta|^2)^{-(i\lambda+n/2)/2}.
\end{eqnarray}

The plan is now to prove an inversion formula and a Plancherel theorem. The following lemma will be very useful in the sequel.

\begin{lemma}\label{mellinlemma}
Let $f \in C_0^{\infty}(\mathbb{R}^n \setminus \{0\})^L$. Then the function $\tilde{f}$ can be written in the form
\begin{eqnarray*}
\tilde{f}(\lambda)=2\pi^{n/2}\frac{\Gamma(-(2i\lambda+(n-2))/2)}{\Gamma(1/2)\Gamma((-2i\lambda+n)/4)}
\mathscr{M}(r \mapsto r^{n/2}f(re_1))(\lambda),
\end{eqnarray*}
where $\mathscr{M}$ is the Mellin transform.
\end{lemma}

\begin{proof}
We start by observing that, since $f$ has compact support
outside the origin, the Mellin transform above admits an entire extension by the Paley-Wiener theorem. 
It thus suffices to prove the statement for $\lambda \in i(n/2-1, \infty)$ by the uniqueness of an analytic continuation. 

By the $L$-invariance of $f$ and of the Lebesgue measure, we have
\begin{eqnarray}
&&\int_{\mathbb{R}^n}f(x)|\langle x,(1,\eta)\rangle|^{-(i\lambda+n/2)}dx .\label{radialisering}\\
&&= \int_{\mathbb{R}^n}f(x)\int_{SO(n)}|\langle gx, (1,\eta) 
\rangle|^{-(i\lambda+n/2)}dgdx.\nonumber
\end{eqnarray}
Consider now the function
\begin{eqnarray}
R(x,y):=\int_{SO(n)}|\langle gx,y \rangle|^{-(i\lambda+n/2)}dg, \qquad x, y \in \mathbb{R}^n.
\end{eqnarray}
It is $SO(n)$-invariant in each variable separately, and it is homogeneous
of degree $-(i\lambda+n/2)$. Hence, 
\begin{eqnarray}
\qquad R(x,y)=|x|^{-(i\lambda+n/2)}|y|^{-(i\lambda+n/2)} \int_{SO(n)}|\langle ge_1,e_1 \rangle|^{-(i\lambda+n/2)}dg. \label{R2}
\end{eqnarray}
The integral on the right hand side can be expressed as an integral over the sphere $S^{n-1}$. Indeed, the fibration
\begin{eqnarray}
p: SO(n) \rightarrow S^{n-1},
p(g)=ge_1
\end{eqnarray}
defines a measure $\sigma$ on $S^{n-1}$ as the pushforward of the normalised Haar measure on $SO(n)$, i.e., $\sigma$ is defined as 
an $SO(n)$-invariant linear functional
on $C(S^{n-1})$ by the equation
\begin{eqnarray}
\int_{S^{n-1}}f(\xi)d\sigma(\xi):=\int_{SO(n)}f(p(g))dg, \qquad f \in C(S^{n-1}). \label{pushforward} 
\end{eqnarray}
By choosing $f$ as a constant function in the above equality, we see that $\sigma$ is the normalised surface measure on $S^{n-1}$.
Applying \eqref{pushforward} to the equality \eqref{R2}, we get
\begin{eqnarray}
R(x,y)=|x|^{-(i\lambda+n/2)}|y|^{-(i\lambda+n/2)}\int_{S^{n-1}}|\zeta_1|^{-(i\lambda+n/2)}d\sigma(\zeta). \label{R-viktig}
\end{eqnarray}
The last integrand depends only on one variable, and hence we can apply the ``Functions of fewer variables''-theorem (cf. \cite{rudinft})
and replace the integral by an integral over the unit interval on the real line.
This yields
\begin{eqnarray*}
&&\int _{S^{n-1}}|\zeta_1|^{-(i\lambda + n/2)} d\sigma(\zeta)\\
&&\frac{2\Gamma(n/2)}{\Gamma(1/2)\Gamma((n-1)/2)}\int_0^1(1-t^2)^{\frac{n-3}{2}}t^{-(i\lambda+n/2)}dt.
\end{eqnarray*}
By performing the change of variables $s=1-t^2$, we obtain
\begin{eqnarray}
&&\frac{2\Gamma(n/2)}{\Gamma(1/2)\Gamma((n-1)/2)}\int_0^1(1-t^2)^{\frac{n-3}{2}}t^{-(i\lambda+n/2)}dt \nonumber\\ 
&=&\frac{\Gamma(n/2)}{\Gamma(1/2)\Gamma((n-1)/2)}\int_0^1s^{\frac{n-1}{2}-1}(1-s)^{-\frac{2i\lambda+(n-2)}{4}-1}ds \nonumber\\
&=&\frac{\Gamma(n/2)}{\Gamma(1/2)\Gamma((n-1)/2)}\beta\left(\frac{n-1}{2},-\frac{2i\lambda+(n-2)}{4}\right) \nonumber\\
&=&\frac{\Gamma(n/2)\Gamma(-(2i\lambda+(n-2))/2)}{\Gamma(1/2)\Gamma((-2i\lambda+n)/4)}, \nonumber
\end{eqnarray}
and hence
\begin{eqnarray}
\int _{S^{n-1}}|\zeta_1|^{-(i\lambda + n/2)} d\sigma(\zeta)=\frac{\Gamma(n/2)\Gamma(-(2i\lambda+(n-2))/2)}{\Gamma(1/2)\Gamma((-2i\lambda+n)/4)}. 
\label{gamma}
\end{eqnarray}
Inserting \eqref{R-viktig} (with $y=(1,\eta)$) and \eqref{gamma} into \eqref{radialisering} gives
\begin{eqnarray*}
&&\int_{\mathbb{R}^n}f(x)|\langle x,(1,\eta)\rangle|^{-(i\lambda+n/2)}dx\\
&=&(1+|\eta|^2)^{-(i\lambda+n/2)/2}\frac{\Gamma(n/2)\Gamma(-(2i\lambda+(n-2))/2)}{\Gamma(1/2)\Gamma((-2i\lambda+n)/4)}\\
&&\times \int_{\mathbb{R}^n}f(x)|x|^{-(i\lambda+n/2)}dx.
\end{eqnarray*}
Finally, we use polar coordinates to compute the integral on the right hand side. Then
\begin{eqnarray*}
\int_{\mathbb{R}^n}f(x)|x|^{-(i\lambda+n/2)}dx=\frac{2\pi^{n/2}}{\Gamma(n/2)}\int_0^{\infty}r^{n/2}f(re_1)r^{-i\lambda}\frac{dr}{r},
\end{eqnarray*}
and hence the lemma is proved.
\end{proof}

\begin{thm}[Inversion formula]
If $f \in C_0^{\infty}(\mathbb{R}^n \setminus \{0\})^L$, then 
\begin{eqnarray*}
f(re_1)=\frac{\Gamma(1/2)}{2\pi^{n/2}}\int_{\mathbb{R}}\tilde{f}(\lambda)
r^{i\lambda-n/2}\frac{\Gamma((-2i\lambda+n)/4)}{\Gamma((-(2i\lambda+(n-2))/2)}d\lambda.
\end{eqnarray*}
\end{thm}

\begin{proof}
By the previous lemma, we can write 
\begin{eqnarray}
\mathscr{M}(r \mapsto r^{n/2}f(re_1))(\lambda):=\tilde{f}(\lambda)b(\lambda).
\end{eqnarray}
By the assumptions on $f$, the inverse Mellin transform is defined for the left hand side since it is in $L^1$, and
the inversion formula for the Mellin transform yields
\begin{eqnarray}
r^{n/2}f(re_1)=\int_{\mathbb{R}}\tilde{f}(\lambda)b(\lambda)r^{i\lambda-1}d\lambda, 
\end{eqnarray}
i.e.,
\begin{eqnarray*}
f(re_1)=
\frac{\Gamma(1/2)}{2\pi^{n/2}}\int_{\mathbb{R}}\tilde{f}(\lambda)
r^{i\lambda-n/2}\frac{\Gamma((-2i\lambda+n)/4)}{\Gamma((-(2i\lambda+(n-2))/2)}d\lambda.
\end{eqnarray*}
\end{proof}

\begin{remark}
Note that this is a somewhat peculiar looking ``Inversion formula''. It does not express the function $f$ as a weighted superposition
of some canonical functions with respect to the Plancherel measure for the given representation. 
This will become clear by the next theorem. The reason that we prove it is rather because
it serves as a means for proving the Plancherel theorem. 
\end{remark}

\begin{thm}[Plancherel theorem] \label{plancherel1}
For all $f \in C_0^{\infty}(\mathbb{R}^n \setminus\{0\})^L$ we have
\begin{eqnarray*}
\int_{\mathbb{R}^n}|f(x)|^2dx=
\int_{\mathbb{R}} |\tilde{f}(\lambda)|^2
\left|\frac{\Gamma(1/2)}{2\pi^{n/2}}
\frac{\Gamma((-2i\lambda+n)/4)}{\Gamma((-(2i\lambda+(n-2))/2)}\right|^2d\lambda.
\end{eqnarray*}
\end{thm}
\begin{proof}
We introduce some temporary notation and write the inversion formula in the simplified form
\begin{eqnarray*}
f(re_1)=\int_{\mathbb{R}}\tilde{f}(\lambda)r^{i\lambda-n/2}\phi(\lambda)d\lambda.
\end{eqnarray*}
By the inversion formula we then have
\begin{eqnarray*}
\int_{\mathbb{R}^n}|f(x)|^2dx&=&\int_{\mathbb{R}^n}f(x)\int_{\mathbb{R}}\overline{\tilde{f}(\lambda)}|x|^{-i\lambda-n/2}
\overline{\phi(\lambda)}d\lambda dx\\
&=&\int_{\mathbb{R}}\overline{\tilde{f}(\lambda)}\int_{\mathbb{R}^n}f(x)|x|^{-i\lambda-n/2}dx
\overline{\phi(\lambda)}d\lambda. 
\end{eqnarray*}
By the proof of Lemma \ref{mellinlemma}, the inner integral can be seen to be equal to 
$\tilde{f}(\lambda)\phi(\lambda)$, and hence
\begin{eqnarray*}
\int_{\mathbb{R}}\overline{\tilde{f}(\lambda)}\int_{\mathbb{R}^n}f(x)|x|^{-i\lambda-n/2}dx
\overline{\phi(\lambda)}d\lambda=\int_{\mathbb{R}}|\tilde{f}(\lambda)|^2|\phi(\lambda)|^2d\lambda,
\end{eqnarray*}
and this concludes the proof.
\end{proof}

\begin{thm} \label{dekomp1}
The operator $T$ extends to a unitary $H$-intertwining operator 
\begin{eqnarray}
U: L^2((\mathbb{R}^n \setminus \{0\})/ \pm1) \rightarrow \int_{\mathbb{R}}\mathcal{H}_{\lambda}d\mu(\lambda),
\end{eqnarray}
where $\mu$ is the measure determined by the identity
\begin{eqnarray*}
\int_{\mathbb{R}}f(\lambda)d\mu(\lambda):=\int_{\mathbb{R}} f(\lambda)
\left|\frac{\Gamma(1/2)}{2\pi^{n/2}}
\frac{\Gamma((-2i\lambda+n)/4)}{\Gamma((-(2i\lambda+(n-2))/2)}\right|^2d\lambda.
\end{eqnarray*}
\end{thm}

\begin{proof}
By Prop. \ref{H-ekv1} and Thm. \ref{plancherel1}, there exists a unique $H$-intertwining extension 
$U: L^2((\mathbb{R}^n \setminus \{0\})/ \pm1) \rightarrow \int_{\mathbb{R}}\mathcal{H}_{\lambda}d\mu(\lambda)$ of $T$. 
The only thing that remains to prove is the surjectivity of $U$.

This follows immediately from the proof of Theorem 9 in \cite{papper1}. Indeed, by the $H$-equivariance of the operator $U$
the action of the commutative Banach algebra $L^1(H)^{\#}$ of left and right $L$-invariant $L^1$-functions on $H$ is 
intertwined. On each subspace $\mathcal{H}_{\lambda}^L$, a function $f \in L^1(H)^{\#}$ acts as a scalar operator, $\hat{f}(\lambda)$.
If we let $v_{\lambda}$ denote the canonical $L$-invariant vector associated with the spherical representation on $\mathcal{H}_{\lambda}$, and
$\omega \in L^2((\mathbb{R}^n \setminus \{0\})/ \pm 1)$ denote an $L$-invariant vector in the minimal $K$-type, then the positive
functional $\Phi$ on $L^1(H)^{\#}$ given by $\Phi(f)=\langle \pi(f)\omega, \omega \rangle$ can be written as the integral
\begin{eqnarray}
\Phi(f)=\int_{\mathbb{R}}\phi_{\lambda}(f)d\mu(\lambda),
\end{eqnarray}
where $\phi_{\lambda}$ is the multiplicative functional $f \mapsto \langle f(\lambda)v_{\lambda}, v_{\lambda} \rangle_{\lambda}$.
The surjectivity now follows from the proof of Theorem 9 in \cite{papper1} by uniqueness of such an integral decomposition of $\Phi$.
\end{proof}

\section{The branching rule: type $I_{nn}$}
We consider the diffeomorphism
$$\Psi: GL(n,\mathbb{C}) \rightarrow \mathbb{C}^* \times SL(n, \mathbb{C})$$ given by 
$$g \stackrel{\Psi}{\mapsto} (\mbox{det}(g),\mbox{det}(g)^{-1/n}g),$$ where we have chosen the branch of the $n$th-root multifunction
determined by the root of unity with the least argument (i.e. in polar coordinates 
$(re^{i\theta})^{1/n}:=r^{1/n}e^{i\theta/n}$). The mapping $\Psi$ has 
inverse 
$$ \Psi^{-1}: (\lambda, h) \mapsto \lambda^{1/n}h.$$ In this case, however, $\Psi$ is not a group homomorphism since the chosen 
branch of the multifunction is not multiplicative. Instead $\Psi$ is multiplicative up to scalar multiples of modulus one. We shall see later 
that we can still use this diffeomorphism to construct representations of $GL(n, \mathbb{C})$ from representations of $SL(n, \mathbb{C})$ and
$\mathbb{C}^*$ respectively.

\subsection{\emph{Some parabolically induced representations of $SL(n, \mathbb{C})$}}

Let $\mathfrak{a}_0=\mathbb{R}e$, where 

\begin{eqnarray*}
e=
\left(\begin{array}{cc}
n-1 & 0 \\
0 & -I_{n-1}
\end{array} \right).
\end{eqnarray*}
Consider the maximal parabolic subalgebra, $\mathfrak{q_0}$, determined by $\mathfrak{a}_0$, with decomposition
$$\mathfrak{q_0}=\overline{\mathfrak{n}_0} \oplus \mathfrak{m}_0 \oplus \mathfrak{a}_0 \oplus \mathfrak{n}_0,$$
where
\begin{eqnarray*}
\mathfrak{n_0}&=&
\left\{ \left(\begin{array}{cccc}
0 & z_1 & \cdots & z_{n-1} \\
0 & 0 & \cdots & 0\\
\vdots & \vdots & \ddots & 0\\
0 & 0 & \cdots & 0
\end{array} \right)
| z_1, \ldots ,z_{n-1} \in \mathbb{C} \right \},
\\
\mathfrak{m_0}&=&
\left\{ \left(\begin{array}{cc}
0 & 0\\
0 & M 
\end{array} \right)
| M \in \mathfrak{sl}(n-1,\mathbb{C}) \right \},
\\
\overline{\mathfrak{n_0}}&=&
\left\{ \left(\begin{array}{cccc}
0 & 0 & \cdots & 0 \\
z_1 & 0 & \cdots & 0\\
\vdots & \vdots & \ddots & 0\\
z_{n-1} & 0 & \cdots & 0
\end{array} \right)
| z_1, \ldots ,z_{n-1} \in \mathbb{C} \right \}.
\end{eqnarray*}
Here the subspace $\mathfrak{m}_0$ is defined by the property
$$Z_{\mathfrak{h'}}(\mathfrak{a}_0)=\mathfrak{a}_0 \oplus \mathfrak{m}_0,$$
and
\begin{eqnarray*}
\mathfrak{n}_0&=&\left \{ X \in \mathfrak{h'} | \, [H,X]=\alpha(H)X, \qquad \forall H \in \mathfrak{a}_0 \right\},\\
\overline{\mathfrak{n}_0}&=&\left \{ X \in \mathfrak{h'} | \, [H,X]=-\alpha(H)X, \qquad \forall H \in \mathfrak{a}_0 \right\}
\end{eqnarray*}
are the generalised root spaces,
where the root $\alpha \in \mathfrak{a_0}^*$ is determined by 
$$\alpha(e)=n.$$
We let $\rho_0$ denote the half sum of the positive roots counted with multiplicity, i.e., 
$$\rho_0=(n-1)\alpha.$$
On the group level we have the corresponding decomposition
$$Q_0=M_0A_0N_0,$$
where
\begin{eqnarray*}
A_0 &=&
\left \{ \left ( \begin{array}{cc}
e^s & 0 \\
0 & qI_{n-1}\\
\end{array} \right) | s,q\in \mathbb{R}, e^sq^{n-1}=1\right\}\\
M_0 &=& 
\left\{ \left(\begin{array}{cc}
1 & 0\\
0 & M 
\end{array} \right)
| M \in GL(n-1, \mathbb{C}) \right \},\\
N_0 &=&
\left\{ \left( \begin{array}{ccccc}
1 & z_1 & z_2 & \cdots & z_{n-1}\\
0 & 1 & 0 & \cdots & 0\\
0 & 0 & 1 & \cdots & 0\\
\vdots & \vdots & 0 & \ddots & 0\\
0 & 0 & 0 & \cdots & 1
\end{array}
\right)| z_1, \ldots ,z_{n-1} \in \mathbb{C} \right\}
\end{eqnarray*}

Consider now the representation $1 \otimes \exp i\lambda \otimes 1$ of the group 
\begin{eqnarray*}
Q_0=M_0A_0N_0.
\end{eqnarray*}
We realise the induced representation 
\begin{eqnarray}
\tau_{\lambda}:=\mbox{Ind}_{Q_0}^{H'}(1 \otimes \exp (i\lambda + \rho_0) \otimes 1)
\end{eqnarray}
in the Hilbert space $L^2(\overline{N_0}, d\overline{n})$.   

The decomposition on the group level $$H'\stackrel{.}{=}\overline{N_0}M_0A_0N_0,$$ gives that 
for $h \in H', \overline{n} \in \overline{N_0}$, $h\overline{n}$ can be factorised as  
\begin{eqnarray*}
&&\left(\begin{array}{cc}
a & b\\
c & d \end{array} \right)
\left( \begin{array}{cc}
1 & 0 \\
z & I_{n-1}\\
\end{array} \right)\\
&&
=\left(\begin{array}{cc}
1 & 0 \\
\frac{c+dz}{a+bz} & I_{n-1} 
\end{array} \right)
\left( \begin{array}{cc}
\frac{a+bz}{|a+bz|} & 0\\
0 & |a+bx|^{1/n-1}(d-\left(\frac{c+dz}{a+bz}\right)b)
\end{array} \right)\\
&&\times
\left(\begin{array}{cc}
|a+bz| & 0\\
0 & |a+bz|^{-1/n-1}I_{n-1} 
\end{array} \right)
\left( \begin{array}{cc}
1 & \frac{b^t}{a+bz}\\
0 & I_{n-1}
\end{array} \right). 
\end{eqnarray*}

Hence, by identifying $L^2(\overline{N_0},d\overline{n})$ with $L^2(\mathbb{C}^{n-1},dm(z))$,
where $dm(z)$ is the Lebesgue measure on $\mathbb{C}^n$,
we obtain the following formula for the action of $H'$ on functions in the representation space:
\begin{eqnarray*}
\tau_{\lambda}(h) f(z)=|az+b|^{-(i\lambda +n)}f \left(\frac{c+dz}{a+bz}\right ),
\end{eqnarray*} 
where $h=\left(\begin{array}{cc}
a & b\\
c & d \end{array} \right)$.

\subsection{\emph{An intertwining operator}}
Recalling the diffeomorphism 
$$\Psi : GL(n, \mathbb{C}) \rightarrow \mathbb{C}^* \times SL(n,\mathbb{C})$$ from the previous section, we can now form 
the representation 
\begin{eqnarray*}
m_{\lambda} \otimes \mbox{Ind}_{Q_0}^{H'}(1 \otimes \exp (i\lambda+\rho_0) \otimes 1),
\end{eqnarray*}
where $m_{\lambda}(c)=|c|^{i\lambda+\rho}$, 
of $\mathbb{C}^* \times SL(n,\mathbb{C})$. This will in fact give a representation of $GL(n, \mathbb{C})$.
Indeed, suppose that we have $g_1, g_2 \in GL(n, \mathbb{C})$. We can write
$$g_1=\lambda_1^{1/n}h_1,$$ and
$$g_2=\lambda_2^{1/n}h_2,$$
with $\lambda_1, \lambda_2 \in \mathbb{C}^*$ and $h_1, h_2 \in SL(n, \mathbb{C})$.
Then $$g_1g_2=(\lambda_1\lambda_2)^{1/n}\xi(\lambda_1,\lambda_2)h_1h_2,$$
where $|\xi(\lambda_1,\lambda_2)|=1$ and hence the mapping
\begin{eqnarray}
\pi_{\lambda}: g \mapsto m_{\lambda} \otimes \mbox{Ind}_{Q_0}^{H'}(1 \otimes \exp (i\lambda+\rho_0) \otimes 1)\circ \Psi(g)
\end{eqnarray}
defines a unitary representation of $GL(n, \mathbb{C})$. We let $\mathcal{H}_{\lambda}$ denote the corresponding representation space.

For $f \in C_0^{\infty}(\mathbb{C}^n)$, we define the function $Tf: \mathbb{C} \times \mathbb{C}^{n-1} \rightarrow \mathbb{C}$ by
\begin{eqnarray*}
Tf(\lambda,\eta):=\int_{\mathbb{C}^n}f(z)|\langle z,(1,\eta)|^{-(i\lambda+n)}dm(z),
\end{eqnarray*}
where the right hand side is to be interpreted using analytic continuation in the variable $\lambda$ 
as in the previous section (eq.\eqref{T-integral} and the following discussion). 
We have the following analog of Prop. \ref{L2}.

\begin{prop} 
For $f \in C_0^{\infty}((\mathbb{C}^n \setminus \{0\}) / \pm 1)$ and $\lambda \in \mathbb{R}$, the functions $T_{\lambda}f$ is
in $L^2(\mathbb{C}^{n-1})$.
\end{prop}
The proof is the same as that of Prop. \ref{L2}.

All the following integral equalities
where the variable $\lambda$ occurs are to be 
thought of as analytic continuations of the corresponding equalities involving convergent integrals. We write $T_{\lambda}f$ for
the function $\eta \mapsto Tf(\lambda,\eta)$.
\begin{prop}
The operator 
\begin{eqnarray}
T_{\lambda}: C_0^{\infty}((\mathbb{C}^n \setminus \{0\})/U(1)) \rightarrow \mathcal{H}_{\lambda}
\end{eqnarray}
is $H$-equivariant.
\end{prop}

\begin{proof}
Take $g \in H$ and write $g=\zeta h$, where $\zeta$ is a diagonal matrix and $h$ has determinant 1. Moreover we write 
$h^{-1}=
\left( \begin{array}{cc}
a & b\\
c & d
\end{array} \right)$.
Then
\begin{eqnarray*}
T_{\lambda}(gf)(\eta)&=&\int_{\mathbb{C}^n}f(g^*z)|\langle(z_1,z'),(1,\eta)\rangle|^{-(i\lambda+n)}|\det g^*|^2 dm(z)\\
&=&\int_{\mathbb{C}^n}f(x)|\langle(z_1,z'),g^{-1}(1,\eta)\rangle|^{-(i\lambda+n)}dm(z)\\
&=&|\zeta|^{i\lambda+n}\int_{\mathbb{C}^n}f(x)|(a+b\eta)z_1+\langle z',c+d\eta \rangle|^{-(i\lambda+n)}dm(z)\\
&=&|\zeta|^{i\lambda+n}\int_{\mathbb{C}^n}f(x)|a+b\eta|^{-(i\lambda+n)}\\
&&\qquad \quad \times |\langle(z_1,z'),(1,(c+d\eta)(a+b\eta)^{-1})\rangle|^
{-(i\lambda+n/2)}dm(z)\\
&=&\pi_{\lambda}(g)T_{\lambda}(f)(\eta).
\end{eqnarray*}
\end{proof}

For an $L$-invariant function $f$, the function $\eta \mapsto Tf(\lambda,\eta)$ is $L$-invariant by the above proposition.
The Cartan-Helgason theorem (\cite{helg2}) therefore allows us to define the function $\tilde{f}$ by
\begin{eqnarray*}
T(\lambda,\eta)=\tilde{f}(\lambda)(1+|\eta|^2)^{-(i\lambda+n)/2}.
\end{eqnarray*}

\begin{lemma} \label{mellinlemma2}
Let $f \in C_0^{\infty}(\mathbb{C}^n \setminus \{0\})^L$. Then the function $\tilde{f}$ can be written in the form
\begin{eqnarray}
\tilde{f}(\lambda)=4\pi^{n}\frac{\Gamma\left(-\frac{i\lambda+(n-2)}{2}\right)}{\Gamma \left(-\frac{i\lambda+(n-2)}{2}+n-1 \right)}
\mathscr{M}(r \mapsto r^{n}f(re_1))(\lambda),
\end{eqnarray}
where $\mathscr{M}$ is the Mellin transform.
\end{lemma}

\begin{proof}
The proof is almost identical to that of Lemma \ref{mellinlemma}. We assume that $\lambda$ is purely imaginary with big enough imaginary 
part. Then
\begin{eqnarray}
&&\int_{\mathbb{C}^n}f(z)|\langle z,(1,\eta)\rangle|^{-(i\lambda+n)}dm(z) \label{radie2} \\
&&= \int_{\mathbb{C}^n}f(z)\int_{SU(n)}|\langle gz, (1,\eta)
\rangle|^{-(i\lambda+n)}dgdm(z). \nonumber
\end{eqnarray}
The inner integral can be written as
\begin{eqnarray}
&&\int_{\mathbb{C}^n}f(z)|\langle z,(1,\eta)\rangle|^{-(i\lambda+n)}dm(z) \label{apa} \\
&&= \int_{\mathbb{C}^n}f(z)\int_{SU(n)}|\langle gz, (1,\eta)
\rangle|^{-(i\lambda+n)}dg \nonumber\\
&&=|z|^{-(i\lambda+n)}(1+|\eta|^2)^{-(i\lambda+n)/2}\int _{S^{2n-1}}|\zeta_1|^{-(i\lambda + n)} d\sigma(\zeta). \nonumber
\end{eqnarray}
The integrand on the right hand side depends only on one variable, and hence we can apply \cite{rudinft}, Prop. 1.4.4. 
This yields
\begin{eqnarray}
&&\int _{S^{2n-1}}|\zeta_1|^{-(i\lambda + n)} d\sigma(\zeta) \label{beta}\\
&&=\frac{n-1}{\pi}\int_U(1-|z|^2)^{n-2}|z|^{-(i\lambda+n)}dm(z), \nonumber
\end{eqnarray}
where $U$ is the unit disc in $\mathbb{C}$. The last integral can be written as
\begin{eqnarray}
&&\frac{n-1}{\pi}\int_U(1-|z|^2)^{n-2}|z|^{-(i\lambda+n)}dm(z) \label{ceta}\\
&=&2\pi(n-1)\int_0^1(1-t)^{(n-1)-1}t^{-(i\lambda+(n-2))/2-1}dt \nonumber\\
&:=&2\pi(n-1)\beta(n-1,-(i\lambda+(n-2))/2) \nonumber \\
&=&2\pi\frac{\Gamma\left(-\frac{i\lambda+(n-2)}{2}\right)}{\Gamma \left(-\frac{i\lambda+(n-2)}{2}+n-1 \right)}. \nonumber
\end{eqnarray}
Using \eqref{apa}, \eqref{beta}, and \eqref{ceta}, the identity \eqref{radie2} can be rewritten in the form
\begin{eqnarray*}
&&\int_{\mathbb{C}^n}f(z)|\langle z,(1,\eta)\rangle|^{-(i\lambda+n)}dm(z)\\
&=&4\pi^n(1+|\eta|^2)^{-(i\lambda+n)/2}
\frac{\Gamma\left(-\frac{i\lambda+(n-2)}{2}\right)}{\Gamma \left(-\frac{i\lambda+(n-2)}{2}+n-1 \right)}\\
&&\times \int_{\mathbb{C}^n}f(z)|z|^{-(i\lambda+n)}dm(z).
\end{eqnarray*}
Using polar coordinates, the integral on the right is given by
\begin{eqnarray*}
\int_{\mathbb{C}^n}f(z)|z|^{-(i\lambda+n)}dm(z)=\frac{2\pi^n}{\Gamma(n)}
\int_0^{\infty}r^nf(re_1)r^{-i\lambda}\frac{dr}{r},
\end{eqnarray*}
and this proves the statement.
\end{proof}
Since $f$ has compact support outside the origin, the right hand side admits an extension to an entire function by the 
Paley-Wiener theorem.

\begin{thm}[Inversion formula]
If $f \in C_0^{\infty}(\mathbb{C}^n \setminus \{0\})^L$, then 
\begin{eqnarray*}
f(re_1)=\frac{1}{4\pi^n}
\int_{\mathbb{R}}\tilde{f}(\lambda)r^{i\lambda-n}\left(-\frac{i\lambda+n-2}{2}\right)_{n-1}d\lambda,
\end{eqnarray*}
where $( \cdot )_k$ denotes the \emph{Pochhammer symbol} defined as
\begin{eqnarray*}
(t)_0&=&1,\\
(t)_k&=&t(t+1)\cdots(t+k-1), \,k \in \mathbb{N}^+.
\end{eqnarray*}
\end{thm}

\begin{proof}
The identity
\begin{eqnarray*}
(x)_k=\frac{\Gamma(x+k)}{\Gamma(x)}
\end{eqnarray*}
shows that the statement of Lemma \ref{mellinlemma2} can be written in the form
\begin{eqnarray}
\mathscr{M}(r \mapsto r^nf(re_1))=\frac{1}{4\pi^n}\left(-\frac{i\lambda+n-2}{2}\right)_{n-1}\tilde{f}(\lambda)
\end{eqnarray}
The inversion formula for the Mellin transform then yields the identity
\begin{eqnarray*}
f(re_1)=\frac{1}{4\pi^n}
\int_{\mathbb{R}}\left(-\frac{i\lambda+n-2}{2}\right)_{n-1}\tilde{f}(\lambda)r^{i\lambda-n}d\lambda.
\end{eqnarray*}
\end{proof}

\begin{thm}[Plancherel theorem]
For all $f \in C_0^{\infty}(\mathbb{C}^n \setminus \{0\})^L$ we have
\begin{eqnarray*}
\int_{\mathbb{C}^n}|f(z)|^2dm(z)=\left(\frac{1}{4\pi^n}\right)^2
\int_{\mathbb{R}} |\tilde{f}(\lambda)|^2
\left|\left(-\frac{i\lambda+n-2}{2}\right)_{n-1}\right|^2d\lambda.
\end{eqnarray*}
\end{thm}
\begin{proof}
For simplicity, we write the inversion formula in the form
\begin{eqnarray*}
f(re_1)=\int_{\mathbb{R}}\tilde{f}(\lambda)r^{i\lambda-n}\phi(\lambda)d\lambda.
\end{eqnarray*}
By the inversion formula we then have
\begin{eqnarray*}
\int_{\mathbb{C}^n}|f(z)|^2dm(z)&=&\int_{\mathbb{C}^n}f(z)\int_{\mathbb{R}}\overline{\tilde{f}(\lambda)}|z|^{-i\lambda-n}
\overline{\phi(\lambda)}dm(z)\\
&=&\int_{\mathbb{R}}\overline{\tilde{f}(\lambda)}\int_{\mathbb{C}^n}f(z)|z|^{-i\lambda-n}dm(z)
\overline{\phi(\lambda)}d\lambda. 
\end{eqnarray*}
By the proof of Lemma \ref{mellinlemma2}, the inner integral can be seen to be equal to 
$\tilde{f}(\lambda)\phi(\lambda)$, and hence
\begin{eqnarray*}
\int_{\mathbb{R}}\overline{\tilde{f}(\lambda)}\int_{\mathbb{C}^n}f(z)|z|^{-i\lambda-n}dm(z)
\overline{\phi(\lambda)}d\lambda=\int_{\mathbb{R}}|\tilde{f}(\lambda)|^2|\phi(\lambda)|^2d\lambda,
\end{eqnarray*}
and this concludes the proof.
\end{proof}

By the same argument that we used to prove Theorem \ref{dekomp1}, we have the following branching law. 
\begin{thm} 
The operator $T$ extends to a unitary $H$-intertwining operator 
\begin{eqnarray}
U: L^2((\mathbb{C}^n \setminus \{0\})/U(1)) \rightarrow \int_{\mathbb{R}}\mathcal{H}_{\lambda}d\mu(\lambda),
\end{eqnarray}
where $\mu$ is the measure determined by the identity
\begin{eqnarray*}
\int_{\mathbb{R}}f(\lambda)d\mu(\lambda):=
\left(\frac{1}{4\pi^n}\right)^2
\int_{\mathbb{R}} f(\lambda)
\left|\left(-\frac{i\lambda+n-2}{2}\right)_{n-1}\right|^2d\lambda.
\end{eqnarray*}
\end{thm}

\begin{remark}
The ideas in this paper could probably be extended to the case of the type $III_n$ bounded symmetric domain
consisting of complex antisymmetric $n \times n$ matrices by realising the corresponding
minimal representation as the Hilbert space \\
$L^2((\mathbb{H}^n \setminus \{0\})/Sp(1))$ and proceeding in an analogous way.
\end{remark}


\begin{thebibliography}{99}


\bibitem{doz}
{Davidson, Mark and {\'O}lafsson, Gestur and Zhang, Genkai},
{\it{Laplace and {S}egal-{B}argmann transforms on {H}ermitian
              symmetric spaces and orthogonal polynomials}},
{J. Funct. Anal.,
 {\bf 204},
 (2003),
1,
157--195}
     

\bibitem{FKsymm}
{Faraut, J. and Kor{\'a}nyi, A.}
{``Analysis on symmetric cones''},
{Oxford Mathematical Monographs,
 Oxford Science Publications,
 The Clarendon Press Oxford University Press,
 New York,
 1994}

\bibitem{helg2}
 {Helgason, Sigurdur}
 {``Groups and geometric analysis''},
 {Mathematical Surveys and Monographs,83, American Mathematical Society,Providence, RI, 2000}
 

\bibitem{knapp1}
{Knapp, Anthony W.}
{``Representation theory of semisimple groups , an overview based on examples''}
{Princeton University Press,
Princeton, NJ,
2001}

\bibitem{Rossi-V}
{Rossi, H. and Vergne, M.}
{{\it Analytic continuation of the holomorphic discrete series of a semi-simple {L}ie group},}
{Acta Math.,
{\bf 136},
(1976),
1-2,
1--59}

\bibitem{rudinft}
{Rudin, Walter}
{``Function theory in the unit ball of {${\bf C}\sp{n}$}''},
{Grundlehren der Mathematischen Wissenschaften [Fundamental
              Principles of Mathematical Science],
241,
Springer-Verlag,
New York,
1980}  


\bibitem{papper1}
{Sepp\"anen, H.}
{{\it Branching of some holomorphic representations of $SO(2,n)$,}}
{Journal of Lie theory (to appear)}


\bibitem{papper3}
{Sepp\"anen, H.}
{{\it Branching laws for minimal holomorphic representations,}}
{Preprint}


\bibitem{dijkpevz}
{van Dijk, G. and Pevzner, M.}
{{\it Berezin kernels of tube domains}},
{J. Funct. Anal.,
{\bf 181},
(2001),
2,
189--208}


\bibitem{wall}
{Wallach, Nolan R.}
{{\it The analytic continuation of the discrete series. {I}, {II}}},
{Trans. Amer. Math. Soc.,
{\bf 251},
(1979),
1--17, 19--37}


\bibitem{zhtams}
{Zhang, Genkai}
{{\it Berezin transform on real bounded symmetric domains}},
{Trans. Amer. Math. Soc.,
{\bf 353},
(2001),
9,
3769--3787 (electronic)}

\bibitem{zhtp}
{Zhang, Genkai}
{{\it Tensor products of minimal holomorphic representations}},
{Represent. Theory,
{\bf 5},
(2001),
164--190 (electronic)}

\bibitem{zhSB}
{Zhang, Genkai}
{{\it Branching coefficients of holomorphic representations and
              {S}egal-{B}argmann transform}},
{J. Funct. Anal.,
{\bf 195},
(2002),
2,
306--349}


\end{thebibliography}
\end{document}